\documentclass{article}
\usepackage{a4wide}
\usepackage[textwidth=125mm,textheight=185mm]{geometry}
\usepackage{color,eucal,enumerate,mathrsfs,amsthm}
\usepackage[normalem]{ulem}
\usepackage{amsmath,amssymb,amsthm,epsfig,bbm}
\usepackage{amsfonts}
\usepackage{dsfont}
\usepackage{mathtools}
\usepackage{leftidx}
\usepackage{graphicx}
\usepackage{esint}
\usepackage{authblk}

\usepackage[pdfborder={0 0 0}]{hyperref}  

\usepackage[latin1]{inputenc}
\usepackage[english]{babel}

\newtheorem{Definition}{Definition}
\newtheorem{Proposition}[Definition]{Proposition}

\newtheorem{Theorem}[Definition]{Theorem}

\theoremstyle{definition}

\allowdisplaybreaks

\renewcommand{\H}{\mathbb{H}}

\newcommand{\R}{\mathbb{R}}

\newcommand{\mm}{{\mbox{\boldmath$m$}}}

\newcommand{\ggamma}{{\mbox{\boldmath$\gamma$}}}

\newcommand{\ppi}{{\mbox{\boldmath$\pi$}}}

\newcommand{\sfd}{{\sf d}}

\newcommand{\sfh}{{\sf h}}

\newcommand{\Kliminf}{K\kern-3pt-\kern-2pt\mathop{\rm lim\,inf}\limits}  % Kuratowski liminf di insiemi
\newcommand{\supp}{\mathop{\rm supp}\nolimits}   % supporto 
\renewcommand{\d}{{\mathrm d}}
\newcommand{\dt}{{\d t}}
\newcommand{\ddt}{{\frac \d\dt}}

\newcommand{\restr}[1]{\lower3pt\hbox{$|_{#1}$}} %restrizione funzione

\newcommand{\la}{\left<}                  % brackets
\newcommand{\ra}{\right>}
\newcommand{\eps}{\varepsilon}  
\newcommand{\nchi}{{\raise.3ex\hbox{$\chi$}}}

\newcommand{\prob}[1]{\mathscr P(#1)}                   % misure di probabilita
\newcommand{\probt}[1]{\mathscr P_2(#1)}                   % misure di probabilita con momento secondo finito
\newcommand{\e}{{\rm{e}}}                           % mappa di valutazione, bisogna mettere `a mano' il tempo t
\renewcommand{\mm}{\mathfrak m}                                %misura di riferimento

\newcommand{\X}{{\rm X}}

\newcommand{\h}{{\sfh}}

\newcommand{\RCD}{{\sf RCD}}
\newcommand{\Ggamma}{{\mathbf\Gamma}}

\newcommand{\HS}{{\lower.3ex\hbox{\scriptsize{\sf HS}}}}
\renewcommand{\H}[1]{{\rm Hess}(#1)}
\renewcommand{\div}{{\rm div}}

\newcommand{\hr}{{\sf r}}
\newcommand{\hR}{{\sf R}}

\title{Second order differentiation formula on $\RCD(K,N)$ spaces}
\author{Nicola Gigli \thanks{SISSA, Via Bonomea 265, 34136 Trieste (Italy). email: ngigli@sissa.it} \quad Luca Tamanini \thanks{Institut f\"ur Angewandte Mathematik, Universit\"at Bonn, Endenicher Allee 60, 53115 Bonn (Germany). email: tamanini@iam.uni-bonn.de}}

\begin{document}

\maketitle

\begin{abstract}
We prove the second order differentiation formula along geodesics in finite-dimensional $\RCD(K,N)$ spaces. Our approach strongly relies on the approximation of $W_2$-geodesics by entropic interpolations and, in order to implement this approximation procedure, on the proof of new (even in the smooth setting) estimates for such interpolations.
\end{abstract}

\tableofcontents

\section{Main result and comments}
This work is about the development of calculus tools in the setting of $\RCD(K,N)$ spaces $(\X,\sfd,\mm)$ with $K\in\R$ and $N\in[1,\infty)$ (see \cite{AmbrosioGigliSavare11-2} for the original definition with $N=\infty$ and \cite{Gigli12} for the case $N < \infty$). The proofs of the announced results are contained in \cite{GigTam18} and, up to technical difficulties, they rely on \cite{GigTam17}, where the same results are obtained for compact $\RCD(K,N)$ spaces.

Recall that an optimal geodesic test plan $\ppi$ on $\X$ is a probability measure on $C([0,1],\X)$ such that $(\e_t)_*\ppi\leq\ C\mm$ for every $t\in[0,1]$ and some $C>0$ and satisfying
\[
\iint_0^1|\dot\gamma_t|^2\,\d t\,\d\ppi(\gamma)=W_2^2\big((\e_0)_*\ppi,(\e_1)_*\ppi\big).
\]
Here $\e_t:C([0,1],\X)\to\X$ is the evaluation map sending $\gamma$ to $\gamma_t$. Any such $\ppi$ is concentrated on constant speed geodesics and for any couple of measures $\mu_0,\mu_1\in\prob\X$ with bounded densities and supports, there is a unique optimal geodesic test plan such that $(\e_0)_*\ppi=\mu_0$, $(\e_1)_*\ppi=\mu_1$.

From the point of view of calculus on metric measure spaces as developed in \cite{AmbrosioGigliSavare11}, the relation between optimal geodesic test plans and standard geodesics is in some sense the same that there is between Sobolev functions and Lipschitz ones. An example of this phenomenon is the following result (a minor variant of a statement in \cite{Gigli13}), which says that we can safely take one derivative of a $W^{1,2}(\X)$ function along an optimal geodesic test plan:

\begin{Theorem}\label{thm:1}
Let $(\X,\sfd,\mm)$ be a $\RCD(K,\infty)$ space, $\ppi$ an optimal geodesic test plan with bounded support (equivalently: such that $\{\gamma_t\,:\, t\in[0,1],\ \gamma\in\supp(\ppi)\}\subset\X$ is bounded) and $h\in W^{1,2}(\X)$. 

Then the map $[0,1]\ni t\mapsto h\circ\e_t\in L^2(\ppi)$ is in $C^1([0,1],L^2(\ppi))$ and we have
\[
\ddt \big( h\circ\e_t\big)=\la\nabla h,\nabla\phi_t\ra\circ\e_t,\\
\]
for every $t\in[0,1]$, where $\phi_t$ is any function such that for some $s\neq t$, $s\in[0,1]$, the function $-(s-t)\phi_t$ is a Kantorovich potential from $(\e_t)_*\ppi$ to $(\e_s)_*\ppi$.
\end{Theorem}

Our main result here is the extension of the above to second order derivatives. Recalling that the second order Sobolev space $H^{2,2}(\X)$ and the corresponding Hessian are defined in \cite{Gigli14}, we have:

\begin{Theorem}\label{thm:main}
Let $(\X,\sfd,\mm)$ be a $\RCD(K,N)$ space, $N<\infty$, $\ppi$ an optimal geodesic test plan with bounded support and $h\in H^{2,2}(\X)$. 

Then the map $[0,1]\ni t\mapsto h\circ\e_t\in L^2(\ppi)$ is in $C^2([0,1],L^2(\ppi))$ and we have
\begin{equation}
\label{eq:th1}
\frac{\d^2}{\d t^2} \big( h\circ\e_t\big)=\H{h}(\nabla\phi_t,\nabla\phi_t)\circ\e_t,\\
\end{equation}
for every $t\in[0,1]$, where $\phi_t$ is as in Theorem \ref{thm:1}. 
\end{Theorem}
Notice that by Theorem \ref{thm:1} we have that such result is really a statement about the $C^1$ regularity of $t\mapsto \la\nabla h,\nabla\phi_t\ra\circ\e_t$. 

Let us collect a couple of equivalent formulations of Theorem \ref{thm:main}. For the first recall that the space of Sobolev vector fields $H^{1,2}_C(T\X)$ as well as the covariant derivative have been defined in \cite{Gigli14}. Then we have:
\begin{Theorem}\label{thm:main2}
Let $(\X,\sfd,\mm)$ be a $\RCD(K,N)$ space, $N<\infty$, $\ppi$ an optimal geodesic test plan with bounded support and $X\in H^{1,2}_C(\X)$. 

Then the map $[0,1]\ni t\mapsto \la X,\nabla \phi_t\ra\circ\e_t\in L^2(\ppi)$ is in $C^1([0,1],L^2(\ppi))$ and we have
\begin{equation}
\label{eq:th2}
\ddt \big(  \la X,\nabla \phi_t\ra\circ\e_t\big)=\nabla X(\nabla\phi_t,\nabla\phi_t)\circ\e_t,
\end{equation}
for every $t\in[0,1]$, where $\phi_t$ is as in Theorem \ref{thm:1}.
\end{Theorem}
From the identity $\nabla(\nabla h)=\H h$ (assuming to identify tangent and cotangent vector fields) we see that Theorem \ref{thm:main2} implies Theorem \ref{thm:main}. For the converse implication notice that Theorem \ref{thm:main} and the Leibniz rule easily provide the correct formula for the derivative of $t\mapsto \la X,\nabla \phi_t\ra\circ\e_t$ for $X=\sum_i\tilde h_i\nabla h_i$, with $(\tilde h_i)\subset L^\infty\cap W^{1,2}(\X)$ and $(h_i)\subset H^{2,2}(\X)$, then conclude by the closure of the covariant derivative.

Another equivalent formulation of Theorem \ref{thm:main}, which is the one we shall actually prove, is:
\begin{Theorem}\label{thm:main3}
Let $(\X,\sfd,\mathfrak{m})$ be a $\RCD(K,N)$ space, $N<\infty$,  $\mu_0,\mu_1 \in \probt \X$ be such that $\mu_0,\mu_1\leq C\mm$ for some $C>0$, with compact supports and let $(\mu_t)$ be the unique $W_2$-geodesic connecting $\mu_0$ to $\mu_1$. Also, let $h\in H^{2,2}(\X)$.

Then the map 
\[
[0,1]\ni \ t\quad \mapsto\quad \int h\,\d\mu_t\ \in\R
\]
belongs to $C^2([0,1])$ and it holds
\begin{equation}
\label{eq:derivate}
\begin{split}
\frac{\d^2}{\d t^2}\int h\,\d\mu_t&=\int \H h(\nabla\phi_t,\nabla\phi_t)\,\d\mu_t,
\end{split}
\end{equation}
for every $t\in[0,1]$, where $\phi_t$ is any function such that for some $s\neq t$, $s\in[0,1]$, the function $-(s-t)\phi_t$ is a Kantorovich potential from $\mu_t$ to $\mu_s$. 
\end{Theorem}
Since for any $W_2$-geodesic as in the statement there is a (unique) optimal geodesic test plan $\ppi$ such that $\mu_t=(\e_t)_*\ppi$ for any $t$, we see that Theorem \ref{thm:main3} follows from Theorem \ref{thm:main} by integration w.r.t.\ $\ppi$. For the converse implication one notices that for any optimal geodesic test plan $\ppi$ with bounded support and $\Gamma\subset C([0,1],\X)$ Borel with $\ppi(\Gamma)>0$, the curve $t\mapsto\ppi(\Gamma)^{-1}(\e_t)_*(\ppi\restr\Gamma)$ fulfils the assumptions of Theorem \ref{thm:main3} with the same $\phi_t$'s as in Theorem \ref{thm:main}. The conclusion then follows by the arbitrariness of $\Gamma$ observing that  $L^2(\ppi)$-derivatives exist for every $t$ if and only if the difference quotients converge in the weak $L^2(\ppi)$-topology for every $t$.

\bigskip

Let us comment about the assumptions in Theorems \ref{thm:main}, \ref{thm:main2}, \ref{thm:main3}:
\begin{itemize}
\item[-] The first order differentiation formula is valid on general $\RCD(K,\infty)$ spaces, while for the second order one we need to assume finite dimensionality. This is due to the strategy of our proof, which among other things uses the Li-Yau inequality.
\item[-] There exist optimal geodesic test plans without bounded support (if $K=0$ or the densities of the initial and final marginals decay sufficiently fast) but in this case the functions $\phi_t$ appearing in the statement(s) are not Lipschitz. As such it seems hard to have $\H{h}(\nabla\phi_t,\nabla\phi_t)\circ\e_t\in L^1(\ppi)$ and thus we can not really hope for anything like \eqref{eq:th1}, \eqref{eq:th2}, \eqref{eq:derivate} to hold: this explains the need of the assumption on bounded supports.
\end{itemize}

Having at disposal the second order differentiation formula is interesting not only at the theoretical level, but also for applications to the study of the geometry of $\RCD$ spaces. For instance, the proofs of both the splitting theorem \cite{Gigli13} and of the `volume cone implies metric cone' \cite{DePhilippisGigli16} in this setting can be greatly simplified by using such formula (in this direction, see \cite{Tamanini17} for comments about the splitting). Also, one aspect of the theory of $\RCD$ spaces which is not yet clear is whether they have constant dimension: for Ricci-limit spaces this is known to be true by a result of Colding-Naber \cite{ColdingNaber12} which uses second order derivatives along geodesics in a crucial way. Thus our result is necessary to replicate Colding-Naber argument in the non-smooth setting (but not sufficient: they also use a calculus with Jacobi fields which as of today does not have  a non-smooth counterpart).

\section{Strategy of the proof}

\subsection{The need of an approximation procedure}

Let us recall that a second order differentiation formula, valid for sufficiently regular curves, has been proved in \cite{Gigli14}:
\begin{Theorem}\label{thm:1i}
Let $(\mu_t)$ be a $W_2$-absolutely continuous curve solving the continuity equation
\begin{equation}\label{eq:cint}
\frac\d{\d t}\mu_t + \div(X_t\mu_t) = 0,
\end{equation}
for some vector fields $(X_t)\subset L^2(T\X)$ in the following sense: for every $f\in W^{1,2}(\X)$ the map $t\mapsto\int f\,\d\mu_t$  is absolutely continuous and it holds
\[
\frac\d{\d t}\int f\,\d\mu_t=\int\la\nabla f,X_t\ra\,\d\mu_t.
\]
Assume that 
\begin{itemize}
\item[(i)] $t \mapsto X_t \in L^2(T{\X})$ is absolutely continuous,
\item[(ii)] $\sup_t\{ \|X_t\|_{L^2} + \|X_t\|_{L^{\infty}}+ \|\nabla X_t\|_{L^2} \} < +\infty$.
\end{itemize} 
Then for $f\in H^{2,2}(\X)$ the map $t\mapsto\int f\,\d\mu_t$ is $C^{1,1}$ and the formula
\begin{equation}
\label{eq:secondsmooth}
\frac{\d^2}{\d t^2}\int f\d\mu_t = \int \H{f}(X_t,X_t) + \la\nabla f,\tfrac\d{\d t} X_t  + \nabla_{X_t} X_t\ra  \d\mu_t
\end{equation}
holds for a.e.\ $t \in [0,1]$.
\end{Theorem}

If the vector fields $X_t$ are of gradient type, so that $X_t=\nabla\phi_t$ for every $t$ and the `acceleration' $a_t$ is defined as
\[
\frac\d{\d t}\phi_t+\frac{|\nabla\phi_t|^2}2=:a_t
\]
then \eqref{eq:secondsmooth} reads as
\begin{equation}
\label{eq:secondsmooth2}
\frac{\d^2}{\d t^2}\int f\d\mu_t = \int \H{f}(\nabla\phi_t,\nabla\phi_t)\,\d\mu_t+\int  \la\nabla f,\nabla a_t\ra  \d\mu_t.
\end{equation}
In the case of geodesics it is well-known that \eqref{eq:cint} holds exactly with $X_t = -\nabla\varphi_t$ for appropriate choices of Kantorovich potentials $\varphi_t$ (see also \cite{GigliHan13} in this direction) and moreover the functions $\varphi_t$ solve (in a sense which we will not make precise here) the Hamilton-Jacobi equation
\begin{equation}
\label{eq:hji}
\frac\d{\d t}\varphi_t = \frac{|\nabla\varphi_t|^2}{2},
\end{equation}
thus in this case the acceleration $a_t$ is identically 0. Hence if the vector fields $(-\nabla\varphi_t)$ satisfied the regularity requirements $(i),(ii)$ in the last theorem, we would easily be able to establish Theorem \ref{thm:main}. However in general this is not the case; informally speaking this has to do with the fact that for solutions of the Hamilton-Jacobi equations we do not have sufficiently strong second order estimates.

In order to establish Theorem \ref{thm:main} it is therefore natural to look for suitable `smooth' approximations of geodesics for which we can apply Theorem \ref{thm:1i} above and then pass to the limit in formula \eqref{eq:secondsmooth}. Given that the source of non-smoothness is in the Hamilton-Jacobi equation it is natural to think at viscous approximation as smoothing procedure: all in all viscous limit is `the' way of approximating the `correct' solution of Hamilton-Jacobi and the Laplacian is well behaved under lower Ricci curvature bounds. However, this does not really work: shortly said, the problem is that not every solution of Hamilton-Jacobi is linked to $W_2$-geodesics, but only those for which shocks do not occur in the time interval $[0,1]$. Since the conclusion of Theorem \ref{thm:main} can only hold along geodesics, we see that we cannot simply use viscous approximation and PDE estimates to conclude (one should incorporate in the estimates the fact that the starting function is $c$-concave, but this seems hard to do).

We shall instead use entropic interpolation, which we now introduce.

\subsection{Entropic interpolation: definition}

Fix two probability measures $\mu_0 = \rho_0\mm$, $\mu_1 = \rho_1\mm$ on $\X$. The Schr\"odinger functional equations are
\begin{equation}
\label{eq:sch10}
\rho_0 = f\,\h_1g \qquad\qquad\qquad\qquad \rho_1 = g\,\h_1f,
\end{equation}
the unknown being the Borel functions $f,g : \X \to [0,\infty)$, where $\h_tf$ is the heat flow starting at $f$ evaluated at time $t$. It turns out that in great generality these equations admit a solution which is unique up to the trivial transformation $(f,g)\mapsto (cf,g/c)$ for some constant $c>0$. Such solution can be found in the following way: let $\hR$ be the measure on $\X^2$ whose density w.r.t.\ $\mm \otimes \mm$ is given by the heat kernel $\hr_t(x,y)$ at time $t=1$ and minimize the Boltzmann-Shannon entropy $H(\ggamma\,|\,\hR)$ among all transport plans $\ggamma$ from $\mu_0$ to $\mu_1$. The Euler equation for the minimizer forces it to be of the form  $f \otimes g\,\hR$ for some Borel functions $f,g : \X \to [0,\infty)$, where $f \otimes g(x,y) := f(x)g(y)$. Then the fact that $f \otimes g\,\hR$ is a transport plan from $\mu_0$ to $\mu_1$ is equivalent to $(f,g)$ solving \eqref{eq:sch10}. 

Once we have found the solution of \eqref{eq:sch10} we can use it in conjunction with the heat flow to interpolate from $\rho_0$ to $\rho_1$ by defining
\[
\rho_t:=\h_tf\,\h_{1-t}g.
\]
This is called {\bf entropic interpolation}. Now we slow down the heat flow: fix $\eps>0$ and by mimicking the above find $f^\eps,g^\eps$ such that
\begin{equation}
\label{eq:sch1}
\rho_0 = f^\eps\,\h_{\eps/2}g^\eps \qquad\qquad \rho_1 = g^\eps\,\h_{\eps/2}f^\eps,
\end{equation}
(the factor $1/2$ plays no special role, but is convenient in computations). Then define 
\[
\rho^\eps_t:=\h_{t\eps/2}f^\eps\,\h_{(1-t)\eps/2}g^\eps.
\]
The remarkable and non-trivial fact here is that as $\eps\downarrow0$ the curves of measures $(\rho^\eps_t\mm)$ converge to the $W_2$-geodesic from $\mu_0$ to $\mu_1$. In order to state our results, it is convenient to introduce the (interpolated) {\bf Schr\"odinger potentials} $\varphi^\eps_t,\psi^\eps_t$ as
\[
\varphi^\eps_t:=\eps\log \h_{t\eps/2}f^\eps\qquad\qquad\qquad\qquad\psi^\eps_t:=\eps\log \h_{(1-t)\eps/2}g^\eps.
\]
In the limit $\eps\downarrow0$ these will converge to forward and backward Kantorovich potentials along the limit geodesic $(\mu_t)$ (see below). In this direction, it is worth to notice that while for $\eps>0$ there is a tight link between potentials and densities, as we trivially have
\[
\varphi^\eps_t+\psi^\eps_t=\eps\log\rho^\eps_t,
\]
in the limit this becomes the well known (weaker) relation that is in place between  forward/backward Kantorovich potentials and measures $(\mu_t)$:
\[
\begin{split}
\varphi_t+\psi_t&=0\qquad\text{on }\supp(\mu_t),\\
\varphi_t+\psi_t&\leq 0\qquad\text{on }\X,
\end{split}
\]
see e.g.\ Remark 7.37 in \cite{Villani09} (paying attention to the different sign convention). By direct computation one can verify that $(\varphi^\eps_t),(\psi^\eps_t)$ solve the Hamilton-Jacobi-Bellman equations
\begin{equation}
\label{eq:hj2}
\frac{\d}{\d t}\varphi^{\varepsilon}_t  =  \frac{1}{2}|\nabla\varphi^{\varepsilon}_t|^2 + \frac{\varepsilon}{2}\Delta\varphi^{\varepsilon}_t\qquad\qquad\qquad\qquad-\frac{\d}{\d t}\psi^{\varepsilon}_t  =  \frac{1}{2}|\nabla\psi^{\varepsilon}_t|^2 + \frac{\varepsilon}{2}\Delta\psi^{\varepsilon}_t,
\end{equation}
thus introducing the functions 
\[
\vartheta^\eps_t:=\frac{\psi^\eps_t-\varphi^\eps_t}2
\]
it is not hard to check that it holds
\begin{equation}
\label{eq:cei}
\frac\d{\d t}\rho^\eps_t+\rm{div}(\nabla\vartheta^\eps_t\,\rho^\eps_t)=0
\end{equation}
and
\[
\frac\d{\d t}\vartheta^\eps_t+\frac{|\nabla\vartheta^\eps_t|^2}2=a^\eps_t,\qquad\qquad\text{where}\qquad a^\eps_t:= -\frac{\varepsilon^2}{8}\Big(2\Delta\log\rho^{\varepsilon}_t + |\nabla\log\rho^{\varepsilon}_t|^2\Big).
\]

\subsection{Entropic interpolations: uniform control and convergence}

With this said, our main results about entropic interpolations can be summarized as follows. Under the assumptions that the metric measure space $(\X,\sfd,\mm)$ is $\RCD(K,N)$, $N<\infty$, and that $\rho_0,\rho_1$ belong to $L^\infty(\X)$ with bounded supports it holds:
\begin{itemize}
\item[-]\noindent\underline{Zeroth order} 
\begin{itemize}
\item[--]\emph{bound} For some $C>0$ we have $\rho^\eps_t\leq C$ for every $\eps\in(0,1)$ and $t\in[0,1]$.
\item[--]\emph{convergence} The curves  $(\rho^\eps_t\mm)$ $W_2$-uniformly converge to the unique $W_2$-geodesic $(\mu_t)$ from $\mu_0$ to $\mu_1$ and setting $\rho_t := \frac{\d\mu_t}{\d\mm}$ it holds $\rho_t^\eps \stackrel{\ast}{\rightharpoonup} \rho_t$ in $L^\infty(\X)$ for all $t \in [0,1]$.
\end{itemize}

\item[-]\noindent\underline{First order}
\begin{itemize}
\item[--]\emph{bound} For any $t \in (0,1]$ the functions $\{\varphi^\eps_t\}_{\eps\in(0,1)}$ are locally equi-Lipschitz. Similarly for the $\psi$'s. 
\item[--]
\emph{convergence} For every sequence $\eps_n\downarrow0$ there is a subsequence - not relabeled - such that for any $t\in(0,1]$ the functions $\varphi^\eps_t$ converge both locally uniformly and in $W^{1,2}_{loc}(\X)$ to a function $\varphi_t$ such that $-t\varphi_t$ is a Kantorovich potential from $\mu_t$ to $\mu_0$. Similarly for the $\psi$'s. 
\end{itemize}

\item[-]\noindent\underline{Second order} For every $\delta\in(0,1/2)$ we have
\begin{itemize}
\item[--]\emph{bound}
\begin{equation}
\label{eq:2i}
\begin{split}
&\sup_{\eps\in(0,1)}\iint_\delta^{1-\delta} \big(|\H{\vartheta^\eps_t}|_\HS^2+\eps^2|\H{\log\rho^\eps_t}|_\HS^2\big)\rho^\eps_t\,\d t\,\d\mm<\infty,\\
&\sup_{\eps\in(0,1)}\iint_\delta^{1-\delta} \big(|\Delta{\vartheta^\eps_t}|^2+\eps^2|\Delta{\log\rho^\eps_t}|^2\big)\rho^\eps_t\,\d t\,\d\mm<\infty.
\end{split}
\end{equation}
Notice that since in general the Laplacian is not the trace of the Hessian, there is no direct link between these two bounds.
\item[--]\emph{convergence} For every function $h\in W^{1,2}(\X)$ with $\Delta h\in L^\infty(\X)$ it holds
\begin{equation}
\label{eq:2ii}
\lim_{\eps\downarrow0}\iint_\delta^{1-\delta}\la\nabla h,\nabla a^\eps_t\ra \rho^\eps_t\,\d t\,\d\mm=0.
\end{equation}
\end{itemize}
\end{itemize}
\bigskip

%As regards the first order convergence statement, we can actually say something more on the limit %functions: setting
%\[
%\varphi_0(x) := \inf_{t \in (0,1]}\varphi_t(x) \qquad\qquad \psi_1(x) := \inf_{t \in [0,1)} \psi_t(x),
%\]
%for all $t \in (0,1]$ it holds $Q_t(-\varphi_0) = -\varphi_t$ and $Q_t(-\psi_1) = -\psi_{1-t}$, where %$Q_t$ denotes the Hopf-Lax semigroup at time $t$.

With the exception of the convergence $\rho^\eps_t\mm \to \mu_t$, all these results are new even on compact smooth manifolds (in fact, even on $\R^d$).

The zeroth and first order bounds are obtained via a combination of Hamilton's gradient estimate and Li-Yau's Laplacian estimate. Similar bounds can also be obtained for the viscous approximation.

The fact that the limit curve $(\mu_t)$ is the $W_2$-geodesic and that the limit potentials are Kantorovich potentials are consequence of the fact that we can pass to the limit in the continuity equation \eqref{eq:cei} and that the limit potentials satisfy the Hamilton-Jacobi equation. Notice that these zeroth and first order convergences are sufficient to pass to the limit in the term with the Hessian in \eqref{eq:secondsmooth2}.

The crucial advantage of dealing with entropic interpolations (which has no counterpart in viscous approximation) is in the second order bounds and convergence results. The key ingredient that allows to obtain these is a formula due to L\'eonard \cite{Leonard13}, who realized that there is a connection between entropic interpolation and lower Ricci bounds; our contribution is the rigorous proof in the $\RCD$ framework of his formal computations:

\begin{Proposition}\label{pro:5}
For any $\eps>0$ the map $t \mapsto H(\mu^{\varepsilon}_t \,|\, \mm)$ belongs to $C([0,1])\cap C^2(0,1)$ and for every $t \in (0,1)$ it holds
\begin{subequations}
\begin{align}
\label{eq:firstder}
\frac{\d}{\d t}H(\mu^{\varepsilon}_t \,|\, \mm) &= \int\la\nabla \rho^\eps_t,\nabla\vartheta^\eps_t\ra\,\d\mm=\frac1{2\eps}\int\big(|\nabla\psi^\eps_t|^2-|\nabla\varphi^\eps_t|^2\big)\rho^\eps_t\,\d\mm,\\
\label{eq:secondder}
\frac{\d^2}{\d t^2}H(\mu^{\varepsilon}_t \,|\, \mm) & =\int\rho^\eps_t\,\d\big(\Ggamma_2(\vartheta^\eps_t)+\tfrac{\eps^2}4\Ggamma_2(\log(\rho^\eps_t))\big)= \frac{1}{2}\int \rho^\eps_t\,\d\big(\Ggamma_2(\varphi^{\varepsilon}_t) + \Ggamma_2(\psi^{\varepsilon}_t)\big).
\end{align}
\end{subequations}
\end{Proposition}

Let us see how to use \eqref{eq:secondder} in the simplified case $K=0$ and $\mm(\X)=1$ to obtain \eqref{eq:2i}. Observe that if $h:[0,1]\to \R^+$ is a convex function, then $-\frac{h(0)}{t}\leq h'(t)\leq \frac{h(1)}{1-t}$ for any $t\in(0,1)$ and thus
\begin{equation}
\label{eq:intderd}
\int_\delta^{1-\delta}h''(t)\,\d t=h'(1-\delta)-h'(\delta)\leq \frac{h(1)}{1-\delta}+\frac{h(0)}{\delta}.
\end{equation}
If $K=0$ we have $\Ggamma_2\geq 0$, so that \eqref{eq:secondder} tells in particular that $t\mapsto H(\mu^{\varepsilon}_t \,|\, \mm) $ is convex for any $\eps>0$, and if $\mm(\X)=1$ such function is non-negative. Therefore \eqref{eq:intderd} gives that for any $\delta\in(0,1/2)$ it holds
\begin{equation}
\label{eq:per2i}
\sup_{\eps\in(0,1)}\int_\delta^{1-\delta}\int \rho^\eps_t\,\d\big(\Ggamma_2(\vartheta^\eps_t)+\tfrac{\eps^2}4\Ggamma_2(\log(\rho^\eps_t))\big)\,\d t\leq \frac{H(\mu_1 \,|\, \mm) }{1-\delta}+\frac{H(\mu_0 \,|\, \mm) }{\delta}<\infty.
\end{equation}
Recalling the Bochner inequalities (\cite{Erbar-Kuwada-Sturm13}, \cite{AmbrosioMondinoSavare13},\cite{Gigli14})
\[
\begin{split}
\Ggamma_2(\eta)&\geq|\H\eta|_\HS^2\mm,\qquad\qquad\qquad\qquad\Ggamma_2(\eta)\geq\frac{(\Delta\eta)^2}N \mm,
\end{split}
\]
we see that \eqref{eq:2i} follows from \eqref{eq:per2i}. Then with some work (see \cite{GigTam17} for the details) starting from \eqref{eq:per2i} we can  deduce \eqref{eq:2ii} which in turn ensures that the term with the acceleration in \eqref{eq:secondsmooth2} vanishes in the limit $\eps\downarrow0$, thus leading to our main result Theorem \ref{thm:main3}.

\bibliographystyle{siam}
{\small
\bibliography{biblio}}

\def\cprime{$'$} \def\cprime{$'$}
\begin{thebibliography}{10}

\bibitem{AmbrosioGigliSavare11}
{\sc L.~Ambrosio, N.~Gigli, and G.~Savar{\'e}}, {\em Calculus and heat flow in
  metric measure spaces and applications to spaces with {R}icci bounds from
  below}, Invent. Math., 195 (2014), pp.~289--391.

\bibitem{AmbrosioGigliSavare11-2}
\leavevmode\vrule height 2pt depth -1.6pt width 23pt, {\em Metric measure
  spaces with {R}iemannian {R}icci curvature bounded from below}, Duke Math.
  J., 163 (2014), pp.~1405--1490.

\bibitem{AmbrosioMondinoSavare13}
{\sc L.~Ambrosio, A.~Mondino, and G.~Savar{\'e}}, {\em Nonlinear diffusion
  equations and curvature conditions in metric measure spaces}.
\newblock Preprint, arXiv:1509.07273, 2015.

\bibitem{ColdingNaber12}
{\sc T.~H. Colding and A.~Naber}, {\em Sharp {H}\"older continuity of tangent
  cones for spaces with a lower {R}icci curvature bound and applications}, Ann.
  of Math. (2), 176 (2012), pp.~1173--1229.

\bibitem{DePhilippisGigli16}
{\sc G.~De~Philippis and N.~Gigli}, {\em From volume cone to metric cone in the
  nonsmooth setting}, Geometric and Functional Analysis, 26 (2016),
  pp.~1526--1587.

\bibitem{Erbar-Kuwada-Sturm13}
{\sc M.~Erbar, K.~Kuwada, and K.-T. Sturm}, {\em On the equivalence of the
  entropic curvature-dimension condition and {B}ochner's inequality on metric
  measure spaces}, Inventiones mathematicae, 201 (2014), pp.~1--79.

\bibitem{Gigli13}
{\sc N.~Gigli}, {\em The splitting theorem in non-smooth context}.
\newblock Preprint, arXiv:1302.5555, 2013.

\bibitem{Gigli14}
\leavevmode\vrule height 2pt depth -1.6pt width 23pt, {\em Nonsmooth
  differential geometry - an approach tailored for spaces with {R}icci
  curvature bounded from below}.
\newblock Accepted at Mem. Amer. Math. Soc., arXiv:1407.0809, 2014.

\bibitem{Gigli12}
\leavevmode\vrule height 2pt depth -1.6pt width 23pt, {\em On the differential
  structure of metric measure spaces and applications}, Mem. Amer. Math. Soc.,
  236 (2015), pp.~vi+91.

\bibitem{GigliHan13}
{\sc N.~Gigli and B.~Han}, {\em The continuity equation on metric measure
  spaces}, Calc. Var. Partial Differential Equations, 53 (2013), pp.~149--177.

\bibitem{GigTam17}
{\sc N.~Gigli and L.~Tamanini}, {\em Second order differentiation formula on
  compact ${RCD}^*({K},{N})$ spaces}.
\newblock Preprint, arXiv:1701.03932, (2017).

\bibitem{GigTam18}
{\sc N.~{Gigli} and L.~{Tamanini}}, {\em Second order differentiation formula
  on ${RCD}^*({K},{N})$ spaces}.
\newblock Preprint, arXiv:1802.02463, (2018).

\bibitem{Leonard13}
{\sc C.~L{\'e}onard}, {\em On the convexity of the entropy along entropic
  interpolations}.
\newblock To appear in ``Measure Theory in Non-Smooth Spaces", Partial
  Differential Equations and Measure Theory. De Gruyter Open, 2017.
  ArXiv:1310.1274v.

\bibitem{Tamanini17}
{\sc L.~Tamanini}, {\em Analysis and geometry of {RCD} spaces via the
  {S}chr{\"o}dinger problem}, PhD thesis, {U}niversit{\'e} {P}aris {N}anterre
  and SISSA, 2017.

\bibitem{Villani09}
{\sc C.~Villani}, {\em Optimal transport. Old and new}, vol.~338 of Grundlehren
  der Mathematischen Wissenschaften, Springer-Verlag, Berlin, 2009.

\end{thebibliography}

\end{document}